\documentclass[10 pt]{amsart}

\usepackage{amsmath}
\usepackage{amsfonts}
\usepackage{amssymb}
\usepackage{amstext}
\usepackage{amsbsy}
\usepackage{amsopn}
\usepackage{amsthm}
\usepackage{amsxtra}
\usepackage{graphicx}
\usepackage{hyperref}
\usepackage{color}
\usepackage[all]{xy}
\usepackage{enumitem}
\usepackage{tikz-cd}
\usepackage{mathtools}
\usepackage{array}
\usepackage{comment}
\usepackage{bm}
\usepackage{soul}

\newtheorem{theorem}{Theorem}[section]
\newtheorem{lemma}[theorem]{Lemma}
\newtheorem{proposition}[theorem]{Proposition}
\newtheorem{corollary}[theorem]{Corollary}

\theoremstyle{definition}
\newtheorem*{definition*}{Definition}

\theoremstyle{remark}
\newtheorem{remark}[theorem]{Remark}
\newtheorem*{remark*}{Remark}

\newcommand{\R}{\mathbb{R}}
\newcommand{\Z}{\mathbb{Z}}

\newcommand{\N}{\mathbb{N}}

\newcommand{\CS}{\mathcal{S}}

\newcommand{\CK}{\mathcal{K}}

\newcommand{\one}{\boldsymbol{1}}

\newcommand{\fin}{\textnormal{fin}}

\newcommand{\diam}{\textnormal{diam}}

\newcommand{\MNS}{\textnormal{M}^{\scriptscriptstyle \textnormal{NS}}}

\newcommand{\inv}[1]{\textnormal{Inv}(\sigma_{n})}
\newcommand{\infinv}[1]{\textnormal{Inv}^{\infty}(\sigma_{n})}

\newcommand{\T}{{\mathbb T}}

\newcommand{\SL}{\textrm{SL}}

\usepackage{todonotes}
\begin{document}

\title{IC-rigidity and generalizations of a theorem of Mahler}
\author{Bryna Kra}
\address{Department of Mathematics, Northwestern University, 2033 Sheridan Road, 
Evanston, IL 60208}
\email{kra@math.northwestern.edu}

\author{Scott Schmieding}
\address{Department of Mathematics, 
107 McAllister Building
University Park, State College, PA 16802}
\email{sks7247@psu.edu}

\thanks{BK was partially supported by the Simons Foundation and NSF grant DMS-2348315 
and SS was partially supported by NSF grant DMS-2247553.}

\begin{abstract}
    Mahler showed that for every irrational $\alpha$ and every base ${b\geq 2}$, there is a single positive integer $n$ such that every digit in $\{0, \dots, b-1\}$ occurs infinitely often in the base-$b$ expansion of $n\alpha$. Furstenberg later observed that dynamical methods can be used to give a different proof.  
    Building on this perspective, we give generalizations of Mahler's theorem in two directions: the multipliers are a syndetic subset of certain multiplicative semigroups, and analogous versions hold in higher dimensions. Our main tool is IC-rigidity, a type of rigidity which is a hybrid of topological and measurable constraints.  
\end{abstract}
\maketitle

\section{Introduction}
For an irrational $\alpha\in\R$, the orbit $\{n\alpha\mod 1: n\geq 1\}$ is dense in the circle $\T = \R/\Z$, and so in particular any prescribed finite block of digits occurs at the beginning of the expansion of some multiple $n\alpha$. Mahler~\cite{Mahler} improved this, showing that for every irrational $\alpha$ and every base $b\geq 2$, there is a single positive integer $n$ such that every digit in $\{0, \dots, b-1\}$ occurs infinitely often in the base-$b$ expansion of $n\alpha$. Using an observation of Furstenberg, Alon and Peres~\cite[Corollary 7.2]{AP} used dynamical methods to give a simple proof of this theorem, and further showed that the set of $n$ with this property can be restricted to certain subsets of $\N$. A quantitative improvement, along with a polynomial version, was given in Berend and Boshernitzan~\cite{Ber-Bosh}.

However, Mahler showed that more holds, meaning that the set of iterates with this property is large in a multiplicative sense. The set $J\subset\N$ is {\em multiplicatively syndetic} if there is a nonempty finite set $F\subset\N$ such that for all $n\in\N$, we have $nF\cap J\neq\emptyset$, where $nF = \{nf: f\in F\}$.  
Mahler showed that the set of $n$ such that $n\alpha$ contains every digit in the base-$b$ expansion is multiplicatively syndetic. 
Equivalently, the density  of these iterates is positive with respect to every multiplicative F\o lner sequence $\Phi$, and the  lower bound can be chosen uniformly over all F\o lner sequences.

We use dynamical methods to give a new proof of Mahler's theorem, and further give a strengthening of conditions on the iterates,  showing that they form a multiplicatively syndetic subset.

\begin{theorem}
\label{th:large-mahler}
Let $1 \le r \in \N$, let $b \ge 2$, and let $\alpha \in \mathbb{R}$ be irrational. Then
$$\{n^{r} : \textrm{ every } 0 \le i \le b-1 \textrm{ appears infinitely often in the base } b \textrm{-expansion of } n^{r}\alpha\}$$
is a multiplicatively syndetic subset of the multiplicative semigroup $\{n^{r} : n \in \mathbb{N}\}$. 
\end{theorem}
It follows immediately that the set of $n$ with this property is infinite. 
Mahler's theorem corresponds to the case $r=1$. 

This theorem is deduced as a corollary of a more general statement in Theorem~\ref{th:mahlergeneral}: the result holds when $\{n^r: n\geq 1\}$ is replaced by any multiplicative semigroup having a largeness property that we call the {\em density Glasner property} (see Section~\ref{sec:density-Mahler}).

The dynamical methods we use build on our previous work in~\cite{KS-IRC}. For a compact metric space $X$, let $\CK(X)$ denote the space of nonempty compact subsets of $X$, endowed with the Hausdorff metric, and we denote this metric by $d(\cdot, \cdot)$.  
Given a countable group (or semigroup) $G$ acting continuously on $X$, an {\em invariant random compact} (IRC) is a Borel probability measure $\mu$ on $\CK(X)$ which is invariant under the induced action of $G$ on $\CK(X)$.  Adding in an assumption that each $g\in G$ acts surjectively on $X$, then the delta mass $\delta_X$ is always an IRC. 

We say the action is \emph{IC-rigid} if the only IRC $\mu$ such that $\mu$-almost every compact is infinite is $\mu = \delta_{X}$. In~\cite{KS-IRC} we showed that the {\em multiplicative action of $\N$ on $\T$}, given by $T_n\colon\T\to\T$ by $T_n(x) = nx\mod 1$ for all $n\geq 1$, has a weaker form of rigidity.  We strengthen this to show that this action is IC-rigid, further showing that the same holds for sparser multiplicative semigroups in $\N$.  

\begin{theorem}
\label{thm:ICrigidityNonT}
For every natural number $r \ge 1$, the action of the multiplicative semigroup $\{n^r:n\in\N\}$ by multiplication on $\T$ is IC-rigid. 
\end{theorem}
In particular, the multiplicative action of $\N$ on $\T$ is IC-rigid. 
As a corollary of this result, we show for every $r \ge 1$, $\varepsilon > 0$, and infinite compact set $Y\subset \T$ that is invariant under multiplication by $p$ for some $p\geq 2$, there are many dilations of the form $T_{n^{r}}(Y)$ which are $\varepsilon$-dense (a precise statement is given in Corollary~\ref{cor:syndetic}). This result, in turn, then allows us to deduce Theorem~\ref{th:large-mahler} in Section~\ref{sec:NonT}. 

Studying the action of ${\rm SL}_d(\Z)$ on $\T^d$, we obtain a higher dimensional version of Mahler's Theorem. Recall for a semigroup $\CS$, a subset $J \subset \CS$ is syndetic if there exists a finite set $F \subset \CS$ such that $\CS=\bigcup_{f \in F}f^{-1}J$ where $f^{-1}J = \{s \in \CS : fs \in J\}$.  To guarantee that  $\delta_{\mathbb{T}^{d}}$ is an invariant IRC, we restrict our attention to the subset $\MNS_{d}(\mathbb{Z})$ of $d \times d$ nonsingular integer matrices. For a $d \times d$ integer matrix $A$ we define $T_{A} \colon \T^{d} \to \T^{d}$ by $T_{A}(x) = Ax$ mod $\mathbb{Z}^{d}$. Given a base $b \ge 2$ and $x \in \mathbb{T}^{d}$, we write
$$x = \sum_{r=1}^{\infty}\frac{c_{r}(x)}{b^{r}}$$
where each $c_{r}(x) = (x_{r}^{(1)},\ldots,x_{r}^{(d)})$ satisfies $0 \le x_{r}^{(i)} \le b-1$ for every $1 \le i \le d$.

\begin{theorem}\label{thm:higherdimmahlerintro}
Let $\alpha \in \mathbb{T}^{d}$ have at least one irrational coordinate. Then the set
$$\{A \in \MNS_{d}(\mathbb{Z}) : \textrm{ for every } w \in \{0,\ldots,b-1\}^{d}\textrm{, } c_{r}(T_{A}(\alpha)) = w \textrm{ for infinitely many } r\}$$
is syndetic in $\MNS_{d}(\mathbb{Z})$.
\end{theorem}
Taking $d=1$, we recover Mahler's theorem. 
Again IC-rigidity is the key ingredient; we prove this holds for certain subgroups of ${\rm SL}_d(\Z)$. 

\begin{theorem}\label{thm:SLICrigid-intro}
Let $\Gamma$ be a subgroup of $\SL_{d}(\mathbb{Z})$ for $d \ge 2$ whose Zariski closure is Zariski-connected, and having the property that the action of $\Gamma$ on $\mathbb{R}^{d}$ is irreducible. Then the action of $\Gamma$ on $\mathbb{T}^{d}$ is IC-rigid.
\end{theorem}

For example, the theorem applies to any subgroup of $\SL_{d}(\mathbb{Z})$ which is Zariski-dense in $\SL_{d}(\mathbb{R})$, 
such as any finite index subgroup of $\SL_d(\Z)$.

Our approach is similar to that used in Theorem~\ref{thm:ICrigidityNonT}, but since the groups in question are not in general amenable we can not make use of F\o{}lner sequences. Instead, we make use of results on random walks for certain subgroups of $\SL_{d}(\mathbb{Z})$ in the work of He and de Saxce~\cite{HeSaxce} and certain bounds from Bulinski and Fish~\cite{Bulinski-Fish}.

\subsection*{Acknowledgements}
We thank Alexander Fish for helpful conversations regarding their results on Glasner properties for subgroups of $SL_{d}(\mathbb{Z})$, and Yair Hartman for helpful conversations regarding a random walk approach for the linear actions.
\subsection*{Methodologies}
AI use was restricted to give a mathematical audit of the paper, resulting in the correction of small typos. 

\section{Multiplicatively syndetic Mahler Theorem}\label{sec:NonT}
\subsection{IC-rigidity for certain multiplicative semigroups of $\N$}
\label{sec:density-Mahler}
We start by proving the action of certain multiplicative semigroups of $\N$ on $\T$ is IC-rigid. For $n \ge 1$, we define $T_{n} \colon \mathbb{T} \to \mathbb{T}$ by setting
$$T_{n}(x) = nx \textrm{ mod } 1.$$
Given $\varepsilon > 0$ we say a set $Y \subset \T$ is $\varepsilon$-dense in $\T$ if for every $x \in \T$ there exists $y \in Y$ such that $d(x,y)<\varepsilon$. We define
$$D_{\varepsilon} = \{Y \in \CK(\mathbb{T}) : Y \textrm{ is } \varepsilon\textrm{-dense in } \mathbb{T}\}.$$
Equivalently, we have that $D_{\varepsilon} = \{Y \in \CK(\mathbb{T}) : d(Y,\mathbb{T})<\varepsilon\}$. We let $\one_{D_\varepsilon}$ denote the indicator function of $D_\varepsilon$. 

\begin{theorem}\label{thm:densityglasnerICrigid}
    If $\CS \subset \mathbb{N}$ is a multiplicative semigroup such that for every infinite 
    $Y \subset \mathbb{T}$ and $\varepsilon > 0$ we have
$$\lim_{N\to\infty}\frac{1}{|\CS \cap [1,N]|}\sum_{n \in \CS \cap [1,N]}\one_{D_{\varepsilon}}(T_nY) = 1,$$
then the action of $\CS$ on $\mathbb{T}$ by multiplication is IC-rigid.
\end{theorem}

  Glasner~\cite{glasner1979} showed that for every infinite set $Y\subset \T$ and $\varepsilon > 0$, some integer dilation $nY$ is $\varepsilon$ dense.  This motivates us to refer to a multiplicative semigroup $\CS\subset\N$ satisfying the hypothesis of Theorem~\ref{thm:densityglasnerICrigid} as having  the {\em density Glasner property}.

\begin{proof}
Suppose $\CS \subset \mathbb{N}$ is a multiplicative semigroup having the density Glasner property. Let $\mu$ be an IRC for $\CS$ acting on $\mathbb{T}$ such that $\mu$-almost every compact set is infinite and let $\varepsilon > 0$. For $N\in\N$ and $\varepsilon > 0$, define  $f_{N,\varepsilon} \colon \CK(\mathbb{T}) \to \mathbb{R}$ by
$$f_{N,\varepsilon}(Y) = \frac{1}{|\CS \cap [1,N]|}\sum_{n \in \CS \cap [1,N]}\one_{D_{\varepsilon}}(T_nY).$$
Since $\mu$ is $T_{n}$-invariant for every $n \in \CS$, we have
$$\mu(D_{\varepsilon}) = \mu(T_{n}^{-1}(D_{\varepsilon})) = \mu\left(\{Y \in \CK(\mathbb{T}) : T_nY \in D_{\varepsilon}\}\right).$$
Then for every $N\in\N$ we have
$$\int_{\CK(\T)} f_{N,\varepsilon} d \mu = \int_{\CK(\T)}\frac{1}{|\CS \cap [1,N]|}\sum_{n \in \CS \cap [1,N]}\one_{D_{\varepsilon}}(T_nY) d\mu = \mu(D_{\varepsilon}).$$
By the assumption on $\CS$ and since $\mu$-almost every $Y$ is infinite, the sequence $f_{N,\varepsilon}$ converges pointwise $\mu$-almost everywhere to the constant function $1$. By the dominated convergence theorem it follows that $\mu(D_{\varepsilon}) = 1$ 
and hence
$$\mu\left(\bigcap_{m=1}^{\infty}D_{\frac{1}{m}}\right) = 1.$$
Since $\bigcap_{m=1}^{\infty}D_{\frac{1}{m}} = \{\T\}$, it follows that $\mu = \delta_{\T}$. 
\end{proof}

We show that some natural multiplicative semigroups in $\N$ have the density Glasner property.  
For $\CS \subset \N$ and $N \ge 1, $ write
$$A_{N}^\CS(t) = \frac{1}{|\CS \cap [1,N]|}\sum_{n \in \CS \cap [1,N]}e^{2 \pi i n t}.$$

\begin{lemma}\label{lemma:fouriercondition}
Suppose $\CS \subset \N$ is a multiplicative semigroup such that the associated $A_N^\CS(t)$ satisfies the following:
\begin{enumerate}
\item
\label{item:lemma1}
If $t \in \mathbb{T}$ is irrational, then $A_{N}^\CS(t) \to 0$ as $N \to \infty$.
\item
\label{item:lemma2}
There exists a function $\kappa\colon\N\to \R^+$ such that $\kappa(R) \to 0$ as $R\to\infty$ and such that for every rational  $\frac{p}{q} \in \mathbb{T}$ written in reduced form with $q > R$, we have $\limsup_{N \to \infty}|A_{N}^\CS(\frac{p}{q})| \le \kappa(R)$.
\end{enumerate}
Then $\CS$ has the density Glasner property.
\end{lemma}
We use the following standard Fourier detection result (see for instance~\cite[Lemma 2.3]{Bulinski-Fish}).
\begin{lemma}
Let $0 < \varepsilon < \frac{1}{2}$. There exists $M= M(\varepsilon) \ge 1$ and $C = C(\varepsilon) > 0$ such that for every $k \ge 1$, if $\{z_{1},\ldots,z_{k}\}$ is not $\varepsilon$-dense in $\T$, then
$$C\sum_{m=1}^{M}\big\vert \frac{1}{k} \sum_{j=1}^{k}e^{2 \pi i m z_{j}}\big\vert^{2}\geq 1.$$
\end{lemma}
Note both $M$ and $C$ depend on $\varepsilon$, but do not depend on either $k$ or on the set $\{z_{1},\ldots,z_{k}\}$.

\begin{proof}[Proof of Lemma~\ref{lemma:fouriercondition}]
Suppose $\CS \subset \mathbb{N}$ is a semigroup satisfying the hypotheses of Lemma~\ref{lemma:fouriercondition}. It suffices to show that for every infinite $Y \subset \mathbb{T}$ and every $\varepsilon > 0$, we have that 
$$\lim_{N\to\infty}\frac{|\{n \in \CS \cap [1,N] : T_nY \textrm{ is not }\varepsilon\textrm{-dense in } \mathbb{T}\}|}{|\CS \cap [1,N]|} = 0.$$
Suppose then that $Y$ is infinite and fix $R$ and $k$. Then we can choose $E_{k} = \{z_{1},\ldots,z_{k}\} \subset Y$ such that $\diam(E_{k}) < \frac{1}{MR}$,  and so $\max_{i,j}\{d(z_{i},z_{j})\} < \frac{1}{MR}$. Then for every $i \ne j$ and every $1 \le m \le M$, we have
$$0 < |m(z_{i}-z_{j})| < \frac{1}{R}.$$
It follows that $m(z_{i}-z_{j})$ is either irrational, or is rational with denominator greater than $R$ (in reduced form). Note
$$\frac{1}{|\CS \cap [1,N]|}\sum_{n \in \CS\cap[1,N]}\lvert\frac{1}{k}\sum_{j=1}^{k}e^{2 \pi i mnz_{j}}\rvert^{2} = \frac{1}{k^{2}}\sum_{i,j=1}^{k}A_{N}^\CS(m(z_{i}-z_{j})).$$
If $T_nY$ is not $\varepsilon$-dense, then $T_nE_{k}$ is not $\varepsilon$-dense and so 
\begin{align*}
    \frac{|\{n \in \CS \cap [1,N]: T_nY \textrm{ is not }\varepsilon\textrm{-dense in } \mathbb{T}\}|}{|\CS \cap [1,N]|} 
   &  \\ \le \frac{C}{|\CS \cap [1,N]|}\sum_{n \in \CS \cap [1,N]} \sum_{m=1}^{M}\big\vert\frac{1}{k}\sum_{j=1}^{k}e^{2 \pi i mnz_{j}}\big\vert^{2}  & =  \frac{C}{k^{2}}\sum_{m=1}^{M}\sum_{i,j=1}^{k}A_{N}^\CS(m(z_{i}-z_{j})).
\end{align*}
Taking absolute values and $\limsup$ as $N \to \infty$ and applying hypotheses~\eqref{item:lemma1} and~\eqref{item:lemma2} to each pair $z_{i}-z_{j}$, we have that 
\begin{align*}
    \limsup_{N \to \infty} & \frac{|\{n \in \CS \cap [1,N] : T_nY \textrm{ is not }\varepsilon\textrm{-dense in } \mathbb{T}\}|}{|\CS \cap [1,N]|}  \\ 
 & \le  CM\frac{1}{k^{2}}(k+k^{2}\kappa(R)) = CM\left(\frac{1}{k}+\kappa(R)\right)
\end{align*}
(note that there are $k$ cases where $i=j$). 
Taking $R \to \infty$ and $k \to \infty$ gives the statement.
\end{proof}

Berend and Peres~\cite{Ber-Per} showed that for every infinite compact set $Y\subset\T$ and every $\varepsilon > 0$, the set of integers $n$ such that $T_nY$ is $\varepsilon$-dense has natural density one. Berend and Boshernitzan~\cite{Ber-Bosh} extended Mahler's theorem to polynomial sequences. 
We extend these results to show that the density conclusion also holds for powers, and this result  immediately implies Theorem~\ref{thm:ICrigidityNonT}.

\begin{theorem}
\label{th:powers}
Let $1 \le r \in \N$ and let $\CS = \{n^{r} : n \ge 1\}$. Then $\CS$ has the density Glasner property, and hence the action of $\CS$ by multiplication on $\T$ is IC-rigid.
\end{theorem}

\begin{proof}
Fix $r \ge 1$ and let $\CS = \{n^{r} : n \ge 1\}$. We check that $\CS$ satisfies Lemma~\ref{lemma:fouriercondition}. For $t$ irrational, by Weyl equidistribution (polynomial form) we have
$$\lim_{N\to\infty}A_{N}(t) = \lim_{N\to\infty}\frac{1}{\lfloor N^{\frac{1}{r}} \rfloor}\sum_{n \le \lfloor N^{\frac{1}{r}} \rfloor}e^{2 \pi i n^{r} t} = 0.$$
Suppose that $t=p/q$ in reduced form. Consider $F_r(p,q) = \frac{1}{q}\sum_{n=1}^{q}e^{2 \pi i p n^{r}/q}$ and let $M_{N} = |\{n : n^{r} \le N\}|$.  Write $M_N = kq+s$ for some $0\leq s < q$.  Then 
$$
A_N(p/q) = \frac{kq}{M_N}F_r(p,q) + \frac{1}{M_N}\sum_{n=1}^s e^{2\pi i pn^r/q}$$
and so 
$$\limsup_{N\to\infty}|A_N(p/q)| = |F_r(p,q)|.$$
By~\cite{HL1,HL2} (also see the discussion in~\cite{LoxtonVaughan}), there exists a constant $C_{r}>0$ (by~\cite{Steckin}, one can take $C_{r} = e^{C(\frac{r}{\varphi(r)})^{2}}$ for some constant $C$) such that
$|qF_r(p,q)| \le C_{r}q^{1-\frac{1}{r}}$ 
and so
$$|F_{r}(p,q)| \le C_{r}q^{-\frac{1}{r}}.$$ 
Using $\kappa(R) = C_{r}R^{-\frac{1}{r}}$, we have $\kappa(R) \to 0$ as $R \to \infty$, and for $q > R$ we have
\begin{equation*}
    \limsup_{N \to \infty}|A_{N}(p/q)| \le \kappa(R).
    \end{equation*}
    It follows by Lemma~\ref{lemma:fouriercondition} that $\CS$ is density Glasner, and hence by Theorem~\ref{thm:densityglasnerICrigid}, the action of $\CS$ is IC-rigid. 
\end{proof}

Applying Theorem~\ref{th:powers} with $r=1$, we deduce the following.
\begin{theorem}\label{th:pldICrigidity}
Suppose $\CS \subset \mathbb{N}$ is a multiplicative semigroup with positive lower density. Then the action of $\CS$ on $\mathbb{T}$ by multiplication is IC-rigid.
\end{theorem}

\begin{proof}
We follow the proof of Theorem~\ref{thm:densityglasnerICrigid}, and so it suffices to show that for $Y$ infinite, we have
$$\lim_{N\to\infty}\frac{|\{n \in \CS \cap [1,N] : T_nY \textrm{ is not }\varepsilon\textrm{-dense in } \mathbb{T}\}|}{|\CS \cap [1,N]|} = 0.$$
Let $Y$ be infinite and $\varepsilon > 0$. Let $B = \{n \in \mathbb{N} : T_nY \textrm{ is not }\varepsilon\textrm{-dense in } \mathbb{T}\}$. Since $\CS$ has positive lower density, there exists $\delta > 0$ such that $|\CS \cap [1,N]| \ge \delta N$ for all sufficiently large $N$. Then
$$\limsup_{N\to\infty}\frac{|B \cap \CS \cap [1,N]|}{|\CS \cap [1,N]|} \le \limsup_{N\to\infty}\frac{|B \cap [1,N]|}{|\CS \cap [1,N]|} \le \limsup_{N\to\infty}\frac{1}{\delta}\frac{|B \cap [1,N]|}{N} = 0, $$
by noting that the limit holds for $\CS=\mathbb{N}$ by the $r=1$ case of Theorem~\ref{th:powers}.
\end{proof}

\subsection{Further corollaries of IC-rigidity for the action of $\N$ on $\T$}
Replacing the use of weak IC-rigidity in~\cite{KS-IRC} by the stronger IC-rigidity given in Theorems~\ref{th:powers} and~\ref{th:pldICrigidity}, we deduce several corollaries.  Though the proofs are analogous, we include them for completeness.

\begin{theorem}\label{thm:invsetsdilations}
Let $Y \subset \mathbb{T}$ be an infinite compact subset invariant under $T_{p}$ for some $p \ge 2$ and let $\CS \subset \N$ be a multiplicative semigroup with the density Glasner property. Then for every F\o{}lner sequence $(F_{m})_{m\in\N}$ of the multiplicative semigroup $\CS$ and every $\varepsilon > 0$, we have
$$\lim_{m\to\infty}\frac{1}{|F_{m}|}\big\vert\{n \in F_{m} : T_nY \textrm{ is } \varepsilon\textrm{-dense in } \mathbb{T}\}\big\vert =1$$
\end{theorem}
\begin{proof}
Let $\mathcal{I}_{p} = \{Z \in \CK(\T) : T_{p}Z = Z \textrm{ and } Z \textrm{ is infinite}\}$. Suppose $Y \in \mathcal{I}_{p}$ for some $p \ge 2$ and set $\mu_{m} = \frac{1}{|F_{m}|}\sum_{n \in F_{m}}\delta_{T_{n}Y}$. 
It suffices to show that the only weak* accumulation point of $\mu_{m}$ is $\delta_{\mathbb{T}}$. Every weak* accumulation point of $\mu_{m}$ is $\CS$-invariant, and by Theorem~\ref{thm:densityglasnerICrigid} the action of $\CS$ on $\mathbb{T}$ is IC-rigid, 
and so it suffices to show that for every weak* accumulation point $\mu$ of the measures $(\mu_{m})_{m\in\N}$, we have $\mu(\CK_{\textrm{fin}}(\mathbb{T})) = 0$.

We show that $\mu(\CK_{\le r}(\T)) = 0$ for every $r \ge 1$. Since $\CK_{\fin}(\T) = \bigcup_{r \ge 1}\CK_{\le r}(\T)$, this  completes the proof.

Fix $r \ge 1$. We claim there exists $\delta_{r} > 0$ such that $d(\mathcal{I}_{p},\CK_{\le r}(\T)) \ge \delta_{r}$. If not, then since $\CK_{\le r}(\T)$ is compact, there exists a finite set $E \in \CK_{\le r}(\T)$ and sequence $Y_{j} \in \mathcal{I}_{p}$ such that $Y_{j} \to E$ in $\CK(\T)$ as $j\to\infty$. 
Since every $Y_{j}\in \mathcal{I}_{p}$ is $T_{p}$-invariant, it follows that $E$ is as well. Choose $k\in\N$ such that $T_{p}^{k}(x) = x$ for every $x \in E$. Then there exists $\delta > 0$ such that $T_{p}^{k}(B(x,\delta)) \cap B(y,\delta) = \emptyset$ for every $x, y \in E$ with $x \ne y$, and $T_{p}^{k}$ acts on $B(x,\delta)$ by a linear expansion by $p^{k}$ for every $x \in E$. Then for some $j\in\N$, we have
$$Y_{j} \subset \bigcup_{x \in E}B(x,\delta).$$
Since $Y_{j}$ is $T_{p}$-invariant, we have $T_{p}^{k}Y_{j} = Y_{j}$. Thus if $z \in Y_{j} \cap B(x,\delta)$ for some $x \in E$, then $T_{p}^{k}(z) \in Y_{j} \cap B(x,\delta)$, and hence $T_{p}^{\ell k}(z) \in Y_{j} \cap B(x,\delta)$ for every $\ell \ge 1$. However since  $T_{p}^{k}$ acts by a linear expansion by $p^{k}$ on $B(x,\delta)$, this is only possible if $z = x$. Thus $Y_{j} \subset E$ contradicting that $Y_{j}$ is infinite and proving the claim.

Since $T_{n}Y \in \mathcal{I}_{p}$ for every $n \in \CS$, it follows that  $d(T_{n}Y,K_{\le r}(\T)) \ge \delta _{r}$ for all such $n$. Letting $B_{\delta_{r}}(\CK_{\le r}(\T)) = \{Y \in \CK(T) : d(Y,\CK_{\le r}(\T)) < \delta_{r}\}$, it follows that $\mu_{m}(B_{\delta_{r}}(\CK_{\le r}(\T))) = 0$ for every $m \ge 1$. Since $B_{\delta_{r}}(\CK_{\le r}(\T))$ is open, this means that $\mu(B_{\delta_{r}}(\CK_{\le r}(\T))) = 0$ and $\mu(\CK_{\le r}(\T)) = 0$.
\end{proof}

\begin{corollary}
\label{cor:syndetic}
Let $Y \subset \mathbb{T}$ be an infinite compact subset invariant under $T_{p}$ for some $p \ge 2$, and let $\CS \subset \N$ be a multiplicative semigroup with the density Glasner property. Then for every $\varepsilon > 0$, the set
$$\{n \in \CS : T_nY \textrm{ is }\varepsilon\textrm{-dense in } \mathbb{T}\}$$
is a multiplicatively syndetic subset of $\CS$.
\end{corollary}

\begin{proof}
It suffices to show that if $J \subset \CS$ satisfies $$\lim_{k\to\infty}\frac{|J \cap F_{k}|}{|F_{k}|} = 1$$ for every F\o{}lner sequence for the multiplicative semigroup $\CS$, then $J$ is multiplicatively syndetic in $\CS$. To see this, if $J$ is not multiplicatively syndetic, then for every finite $F\subset \CS$, there exists some $r\in \CS$ such that 
$rF\subset \CS \setminus J$. Let $F_m$ be a F\o{}lner sequence for the multiplicative semigroup $\CS$. For each $m$, pick $r_m\in \CS$ such that $F_m r_m \subset \CS \setminus J$ and set 
$\Psi_m = F_m r_m$.  Then $(\Psi_m)$ is also a F\o{}lner sequence for $\CS$, and by construction lies entirely in $\mathbb{N} \setminus J$.
\end{proof}

\subsection{Generalizing Mahler's Theorem}
\begin{theorem}\label{th:mahlergeneral}
Let $\CS \subset \N$ be a multiplicative semigroup having the density Glasner property, let $1 \le b \in \N$, and let $\alpha \in \mathbb{R}$ be irrational. Then
$$\{n \in \CS : \textrm{ every } 0 \le i \le b-1 \textrm{ appears infinitely often in the base } b \textrm{-expansion of } n\alpha\}$$
is a multiplicatively syndetic subset of $\CS$.
\end{theorem}

\begin{proof}
Let $\alpha \in \mathbb{R}$ be irrational and $2 \le b \in \mathbb{N}$. Let $\mathcal{B}$ be the set of limit points of the set $\{b^{k}\alpha \textrm{ mod 1} : k \ge 1\} \subset \mathbb{T}$. The set $\mathcal{B}$ is $T_{b}$-invariant. We claim $\mathcal{B}$ is infinite. Indeed, suppose instead that $\mathcal{B}$ were finite. Then since $\mathcal{B}$ is $T_{b}$-invariant, there exists $k\geq 1$ such that $T_{b}^{k}$ pointwise fixes each element of $\mathcal{B}$, and every element of $\mathcal{B}$ is rational. Choosing small neighborhoods of each point in $\mathcal{B}$ and using the fact that $T_{b}$ locally expands, it follows that $T_{b}^{j}(\alpha) \in \mathcal{B}$ for some $j$, which is a contradiction since $\alpha$ is irrational.

Set $\varepsilon = \frac{1}{2b}$ and define
$$J = \{n \in \CS : T_{n}\mathcal{B} \textrm{ is } \varepsilon\textrm{-dense in } \mathbb{T}\}.$$
If $m \in J$ then $m\alpha$ in base $b$ contains every $0 \le i \le b-1$ infinitely often. Indeed, if $m \in J$, then for every $0 \le i \le b-1$ we have $T_{m}\mathcal{B} \cap (\frac{i}{b},\frac{i+1}{b}) \ne \emptyset$. Since each point in $T_{m}\mathcal{B}$ is an accumulation point of $\{b^{n}m\alpha \textrm{ mod } 1 : n \ge 1\}$, it follows that infinitely many points from $\{b^{n}m\alpha \textrm{ mod } 1: n \ge 1\}$ belong to $(\frac{i}{b},\frac{i+1}{b})$. Thus the result follows from Corollary~\ref{cor:syndetic}. 
\end{proof}

We refer to the proofs of~\cite[Theorem 5.7, Corollary 5.8]{KS-IRC} for the proofs of the next two. 

\begin{corollary}
\label{cor:dilate-revised}
Let $Y \subset \mathbb{T}$ be an infinite compact subset invariant under $T_{p}$ for some $p \ge 2$, and suppose $\CS \subset \N$ is a multiplicative semigroup with the density Glasner property. Then for every F\o{}lner sequence $(F_{m})_{m\in\N}$ of $\CS$, there exists a set $J \subset \CS$ such that
$$\lim_{m\to\infty}\frac{|J \cap F_{m}|}{|F_{m}|} = 1$$
and
$$\lim_{n \to \infty, n \in J}T_nY = \mathbb{T} \quad \textrm{ in } \quad \CK(\mathbb{T}).$$
\end{corollary}

\begin{corollary}
\label{cor:mult-synd}
Let $Y \subset \mathbb{T}$ be an infinite compact subset invariant under $T_{p}$ for some $p \ge 2$ and suppose $\CS \subset \N$ is a multiplicative semigroup with the density Glasner property. Let $W \subset \CS$ be finite and fix $\varepsilon > 0$. Then for every F\o{}lner sequence $(F_{m})_{m\in\N}$ of the multiplicative semigroup $\CS$, we have 
$$\lim_{m\to\infty}\frac{1}{|F_{m}|}\lvert \{n \in F_{m} : T_{wn}Y \textrm{ is }\varepsilon \textrm{-dense in }\mathbb{T} \textrm{ for all } w \in W\}| = 1.$$
\end{corollary}

\begin{remark}
Using Corollary~\ref{cor:mult-synd} with $\CS = \mathbb{N}$ we also have the following slightly stronger version of Theorem~\ref{th:mahlergeneral}: Given $b,k \ge 2$ and $\alpha \in \mathbb{R}$ irrational, the set of $n \in \mathbb{N}$ such that every $0 \le i \le b-1$ appears infinitely often in the base $b$-expansions of each of $n\alpha, 2n\alpha, \ldots, kn\alpha$ is multiplicatively syndetic.
\end{remark}

We also have the following.
\begin{corollary}
Let $Y \subset \mathbb{T}$ be an infinite compact subset invariant under $T_{p}$ for some $p \ge 2$ and let $\varepsilon > 0$. There exists a finite set $F \subset \mathbb{N}$ such that for every $n \ge 1$, there exists $f \in F$ such that the set 
$$J(n,\varepsilon) = \{k \in \mathbb{N} : T_{fn^{k}}Y \textrm{ is }\varepsilon\textrm{-dense in } \mathbb{T}\}$$
has positive upper density.
\end{corollary}

\begin{proof}
This is a consequence of Corollary~\ref{cor:syndetic} with $\mathcal{S} = \mathbb{N}$. Indeed, given such $Y$ and $\varepsilon > 0$, by Corollary~\ref{cor:syndetic} the set
$$J = \{m \in \mathbb{N} : T_mY \textrm{ is }\varepsilon\textrm{-dense in } \mathbb{T}\}$$
is multiplicatively syndetic. Choose $F=\{f_{1},\ldots,f_{k}\}$ such that for every $m \in \mathbb{N}$ we have $mF \cap J \ne \emptyset$. Fix $n \ge 1$. Then for every $j \ge 1$, there exists $1 \le i \le k$ such that $f_{i}n^{j} \in J$, and so
$$\bigcup_{i=1}^{k}\{j : f_{i}n^{j} \in J\} = \mathbb{N}.$$
It follows that there exists some $f \in F$ such that
\begin{equation*}
    \overline{d}\left(\{j : fn^{j} \in J\}\right) \ge \frac{1}{|F|}>0, 
    \end{equation*}
where $\overline{d}$ denotes the upper density. 
\end{proof}

\subsection{Further IC-rigid actions}
We remark that there are other mechanisms for producing IC-rigid actions, and note one of them here. 
Recall an action of $G$ on $X$ is {\em extremely proximal} if for every $X \ne Y \in \CK(X)$ there exists $x \in X$ such that $\{x\} \in \overline{\{gY : g \in G\}}$.  See Glasner~\cite{glasner1974} for further background. 
\begin{proposition}
If $(X,G)$ is extremely proximal and minimal with $X$ infinite, then it is IC-rigid.
\end{proposition}

\begin{proof}
Suppose $(X,G)$ is extremely proximal and also minimal, and let $X \ne Y \in \CK(X)$. Minimality together with extreme proximality implies for every $x \in X$ we have $\{x\} \in \overline{\{gY : g \in G\}}$. It follows that if $\emptyset \ne U \subset X$ is open then there exists $g \in G$ such that $gY \subset U$.

Fix $Z \in \CK(X)$ with $Z \ne X$ and write $E(Z) = \{Y \in \CK(X) : Y \subset Z\}$. Given $N \ge 1$, since $X$ is infinite, we may choose nonempty disjoint open sets $U_{1},\ldots,U_{N}$ in $X$. By the above, for each $1 \le i \le N$ there exists $g_{i} \in G$ such that $g_{i}Z \subset U_{i}$, and hence the sets $g_{i}Z$ are pairwise disjoint. Suppose $\mu$ is an IRC for $(X,G)$. Using invariance of $\mu$ we have
$$N \mu(E(Z)) = \sum_{i=1}^{N}\mu(E(Z)) \le \sum_{i=1}^{N}\mu(g_{i}E(Z))  \le 1$$
since each of the $g_{i}E(Z)$ are disjoint. Since this holds for all $N \ge 1$ we get that $\mu(E(Z)) = 0$. If $V_{i}$ is a countable basis for $X$, then
$$\CK(X) \setminus \{X\} \subset \bigcup_{i=1}^{\infty}E(X \setminus V_{i})$$
and it follows that $\mu=\delta_{X}$.
\end{proof}

\section{Subgroups of $\SL_{d}(\mathbb{Z})$}\label{sec:SLdZ}
\subsection{IC-rigidity for certain subgroups of $SL_{d}(\mathbb{Z})$}
Throughout, we fix some $2 \le d \in \mathbb{N}$.  
We slightly abuse notation and write 
$$D_{\varepsilon} = \{Y \in \CK(\mathbb{T}^{d}) : d(Y,\mathbb{T}^{d}) < \varepsilon\}.$$
\begin{theorem}\label{thm:SLICrigid}
Let $\Gamma$ be a subgroup of $\SL_{d}(\mathbb{Z})$ whose Zariski closure is Zariski-connected, and having the property that the action of $\Gamma$ on $\mathbb{R}^{d}$ is irreducible. Then the action of $\Gamma$ on $\mathbb{T}^{d}$ is IC-rigid.
\end{theorem}
 We make a slight abuse of notation and write the action of $\Gamma$ on $\T^d$ as $gE$ for $g\in \Gamma$ and $E\subset \T^d$. 

For the remainder of this section, we fix a subgroup $\Gamma$ of $\SL_{d}(\mathbb{Z})$ whose Zariski closure is Zariski-connected, and having the property that the action of $\Gamma$ on $\mathbb{R}^{d}$ is irreducible. We let $\rho$ be a fully-supported probability measure on $\Gamma$ satisfying $\sum_{g \in \Gamma}\|g\|^{\varepsilon}\rho(g) < \infty$ for some $\varepsilon > 0$, and write $\rho^{*n}$ for the $n$-fold convolution.

Before stating the main technical result, we recall a result that is implicit in the proof of~\cite[Theorem A]{Bulinski-Fish}. We include the proof for convenience. 

\begin{lemma}\label{lemma:BFbounds}
Let $\varepsilon > 0$. There exists $C_{\varepsilon} < \infty$ and $\eta > 0$ such that for every $E = \{x_{1},\ldots,x_{k}\} \subset \mathbb{T}^{d}$, we have
$$\limsup_{n \to \infty}\rho^{*n}\left(\{g \in \Gamma : gE \textrm{ is not } \varepsilon\textrm{-dense in } \mathbb{T}^{d}\}\right) \le C_{\varepsilon}(k^{-\eta} + k^{-1}),$$
and the constant $C_\varepsilon$ does not depend on $k$. 
\end{lemma}
\begin{proof}
Set $M = \lceil d/\varepsilon \rceil$ and define
$$\Phi_{E}(g) = \sum_{v \in B(M)}\sum_{i,j=1}^{k}e^{2 \pi i (v \cdot g(x_{i}-x_{j}))}$$
where $B(M) = \{0 \ne v \in \mathbb{Z}^{d} : \|v\|_{\infty} \le M\}$. Note $\Phi_{E}(g)$ is real and nonnegative, since it is also equal to $\sum_{v \in B(M)}\lvert\sum_{i=1}^{k}e^{2 \pi i v \cdot g(x_{i})}\rvert^{2}$. Using~\cite[Lemma 2.3]{Bulinski-Fish}, we have that if $gE$ is not $\varepsilon$-dense, then
\begin{equation}\label{eqn:notdense}
k^{2} \le C_{d}\varepsilon^{-d}\Phi_{E}(g)
\end{equation}
for some constant $C_{d}$. We consider two cases. If $x_{i}-x_{j}$ is irrational, by the assumptions on $\Gamma$ we can apply~\cite[Theorem 1.2]{HeSaxce} (see the reformulation in~\cite[Lemma 3.2]{Bulinski-Fish}) 
and we have that 
$$\lim_{n\to\infty}\int_{\Gamma}e^{2 \pi i (v \cdot g(x_{i}-x_{j}))}\, d \rho^{*n} = 0.$$
On the other hand, if $x_{i}-x_{j}$ is rational, it has finite order in $\mathbb{T}^{d}$, say order $q$. Again using the assumptions on $\Gamma$ we can apply~\cite[Lemma 3.2]{Bulinski-Fish} to get $c > 0$ (which does not depend on $x_{i}-x_{j}$) such that for each $v$ we have
\begin{equation}
    \label{eq:bound-integral}
\limsup_{n \to \infty} \int_{\Gamma} \lvert e^{2 \pi i (v \cdot g(x_{i}-x_{j}))} \rvert \, d\rho^{*n} \le 2\|v\|q^{-\frac{1}{c}}.
\end{equation}
Without loss of generality, we assume $c \ge 1$. Let $h_{q} = \vert\{(i,j) : x_{i}-x_{j} \textrm{ has order } q\}\vert$. Then~\cite[Proposition 2.1]{Bulinski-Fish} gives
\begin{equation}
    \label{eq:bound-sum}
    \sum_{q \ge 2}h_{q}q^{-\beta} \le Ck^{2-\frac{\beta}{d+1}}.
    \end{equation}
Let $\beta = 1/c$ and set $\eta = \frac{\beta}{d+1}$. 
Each of the  $k$-many $i=j$ cases in the definition of $\Phi_{E}(g)$ adds 1 to the sum for each $v$, in total contributing $|B(M)|k$. Since $|B(M)|$ is finite, applying the  bounds in~\eqref{eq:bound-integral} and~\eqref{eq:bound-sum} we obtain
\begin{equation}\label{eqn:intbound}
\limsup_{n \to \infty}\int_{\Gamma} \Phi_{E}(g)\, d \rho^{*n} \le C^{\prime}_{\varepsilon}\left(k^{2-\eta}+k\right)
\end{equation}
for some $C^{\prime}_{\varepsilon}>0$. By Markov's inequality, it follows that 
$$\rho^{*n}\left( \{g : \Phi_{E}(g) \ge C_{d}^{-1}\varepsilon^{d}k^{2}\}\right) \le \frac{C_{d}\varepsilon^{-d}}{k^{2}}\int_{\Gamma} \Phi_{E}(g)\, d\rho^{*n}.$$
Combining this with~\eqref{eqn:notdense} and~\eqref{eqn:intbound}, we have
\begin{align*}\limsup_{n \to \infty}\rho^{*n} & \left( \{g :  gE \textrm{ not } \varepsilon\textrm{-dense}\}\right) 
\le   \limsup_{n \to \infty} \rho^{*n}\left( \{g : k^{2} \le C_{d}\varepsilon^{-d}\Phi_{E}(g)\} \right)  \\ 
\le & \limsup_{n \to \infty}\frac{C_{d}\varepsilon^{-d}}{k^{2}}\int_{\Gamma} \Phi_{E}(g)\, d\rho^{*n}  \\ 
\le  & \frac{C_{d}\varepsilon^{-d}}{k^{2}}C_{\varepsilon}^{\prime}\left(k^{2-\eta}+k\right) =C_{\varepsilon}\left(k^{-\eta}+k^{-1}\right).  \qedhere
\end{align*}
\end{proof}

This allows us to prove the main technical ingredient for showing IC-rigidity of the action of $\Gamma$ on $\T^d$. 

\begin{proposition}\label{prop:randomepsilondense}
If $d\geq 2$, $Y \subset \mathbb{T}^{d}$ is compact and infinite,  and  $\varepsilon > 0$, then we have 
$$\lim_{n\to\infty}\rho^{*n}\left(\{g \in \Gamma : gY \textrm{ is } \varepsilon\textrm{-dense in } \mathbb{T}^{d}\}\right) = 1.$$
\end{proposition}
\begin{proof}
Given infinite compact $Y \subset \mathbb{T}^{d}$ and $\varepsilon > 0$, for each $k\geq 1$ we may select $F_{k}=\{x_{1},\ldots,x_{k}\} \subset Y$. If $gY$ is not $\varepsilon$-dense, then $gF_{k}$ is not $\varepsilon$-dense. By Lemma~\ref{lemma:BFbounds}, it follows that for each $k\geq 1$ we have that 
$$\limsup_{n \to \infty}\rho^{*n}\left(\{g \in \Gamma : gY \textrm{ is not } \varepsilon\textrm{-dense in } \mathbb{T}^{d}\}\right) \le C_{\varepsilon}(k^{-\eta} + k^{-1}).$$
Taking $k \to \infty$, we obtain
$$\limsup_{n \to \infty}\rho^{*n}\left(\{g \in \Gamma : gY \textrm{ is not } \varepsilon\textrm{-dense in } \mathbb{T}^{d}\}\right) = 0$$
and hence  
\begin{equation*}
   \lim_{n\to\infty} \rho^{*n}\left(\{g \in \Gamma : gY \textrm{ is } \varepsilon\textrm{-dense in } \mathbb{T}^{d}\}\right) = 1. \qedhere\end{equation*}
\end{proof}
We use this to complete the proof of  Theorem~\ref{thm:SLICrigid}.

\begin{proof}[Proof of Theorem~\ref{thm:SLICrigid}]
Suppose $\mu$ is an IRC for the action of $\Gamma$ on $\mathbb{T}^{d}$ such that $\mu$-almost every compact is infinite and let $\varepsilon > 0$. Recall $D_{\varepsilon} = \{Y \in \mathcal{K}(\mathbb{T}^{d}) : Y \textrm{ is }\varepsilon\textrm{-dense in } \mathbb{T}^{d}\}$.
For $n \ge 1$, define $f_{n,\varepsilon} \colon \mathcal{K}(\mathbb{T}^{d}) \to \mathbb{R}$
by
$$f_{n,\varepsilon}(Y) = \rho^{*n}\left(\{g \in \Gamma : gY \in D_{\varepsilon}\}\right).$$
Then 
\begin{align*}
    \int_{\mathcal{K}(\mathbb{T}^{d})}f_{n,\varepsilon}(Y)\, d\mu(Y) = & \int_{\mathcal{K}(\mathbb{T}^{d})} \rho^{*n}\left( \{g : gY \in D_{\varepsilon}\} \right)\, d\mu(Y)\\ = & \int_{\mathcal{K}(\mathbb{T}^{d})} \int_{\Gamma}\mathbf{1}_{D_{\varepsilon}}(gY)\, d\rho^{*n} \, d\mu(Y),
\end{align*}
and by Fubini's theorem this is equal to
\begin{align*}
\int_{\Gamma} \int_{\mathcal{K}(\mathbb{T}^{d})}\mathbf{1}_{D_{\varepsilon}}(gY)\, d\mu \, d\rho^{*n} = & \int_{\Gamma} \mu\left(\{Y : gY \in D_{\varepsilon}\}\right) \, d\rho^{*n} \\ = & \int_{\Gamma} \mu(g^{-1}D_{\varepsilon})\, d\rho^{*n}.
\end{align*}
Since $\mu$ is $\Gamma$-invariant, this last term is equal to
$$\int_{\Gamma} \mu(D_{\varepsilon})\, d\rho^{*n} = \mu(D_{\varepsilon}).$$
By assumption, $\mu$-almost every $Y\in\CK(\T^d)$ is infinite, and so by Proposition~\ref{prop:randomepsilondense}, $f_{n,\varepsilon}(Y)$ converges to $1$ for $\mu$-almost every $Y\in\CK(\T^d)$. It follows by the dominated convergence theorem that $\mu(D_{\varepsilon}) = 1$. Thus
$$\mu\left(\bigcap_{m=1}^{\infty}D_{\frac{1}{m}}\right) = 1.$$
Since $\bigcap_{m=1}^{\infty}D_{\frac{1}{m}} = \{\mathbb{T}^{d}\}$, we have $\mu=\delta_{\mathbb{T}^{d}}$.
\end{proof}

\subsection{Corollaries
and a higher-dimensional Mahler Theorem}
Fix $d \ge 2$. Given $p \ge 2$, we (in a minor abuse of notation) write $T_{p} \colon \mathbb{T}^{d} \to \mathbb{T}^{d}$ where $T_{p}(x)=px$ (since $d$ is fixed throughout, there is no ambiguity as to the space on which $T_{p}$ acts). Define 
$$\mathcal{I}_{p} = \{Y \in \mathcal{K}(\mathbb{T}^{d}) : Y \textrm{ is infinite and }T_{p}Y=Y\}.$$
\begin{lemma}
Both of the following hold.
\begin{enumerate}
\item\label{item:closed}
The set $\mathcal{I}_{p}$ is closed in $\CK(\mathbb{T}^{d})$ and hence compact. 
\item\label{item:large}
For every $r \ge 1$, there exists $\delta_{r} > 0$ such that $d(\mathcal{I}_{p},\mathcal{K}_{\le r}(\mathbb{T}^{d})) \ge \delta_{r}$.
\end{enumerate}
\end{lemma}
\begin{proof}
To  check part~\eqref{item:closed}, suppose there is a sequence  $Y_{j} \in \mathcal{I}_{p}$ such that $Y_j\to Y$ and $Y$ is  finite. Then $T_{p}Y=Y$ and we may choose $k\geq 1$ such that $T_{p}^{k}$ acts by the identity on $Y$, and choose $\delta > 0$ such that $T_{p}^{k}(B_{\delta}(x)) \cap B_{\delta}(y) = \emptyset $ for all $x \ne y \in Y$.  It follows that  $Y_{j} \subset \bigcup_{x \in Y}B_{\delta}(x)$ for $j$ sufficiently large. But then by the invariance $T_{p}^{k}Y_{j} = Y_{j}$, if $z \in Y_{j} \cap B_{\delta}(x)$, then the choice of $\delta$ implies $T_{p}^{k \ell}(z) \in Y_{j} \cap B_{\delta}(x)$ for all $\ell \ge 0$.  Then the local expansion of $T_{p}$ forces $z=x$,  and so $Y_{j} \subset Y$, a contradiction. 

The proof of part~\eqref{item:large} is the same as used in the proof of Theorem~\ref{thm:invsetsdilations}.
\end{proof}

We use this to obtain the analog of Theorem~\ref{th:mahlergeneral} in the higher-dimensional setting. To guarantee that  $\delta_{\mathbb{T}^{d}}$ is an IRC, we restrict our attention to the subset $\MNS_{d}(\mathbb{Z})$ of $d \times d$ nonsingular integer matrices.

\begin{corollary}\label{cor:matricessyndetic}
Suppose $Y \in \mathcal{I}_{p}$. Then for every $\varepsilon > 0$, the set
$$J_{Y} = \{A \in \MNS_{d}(\mathbb{Z}) : T_{A}Y \textrm{ is } \varepsilon\textrm{-dense in } \mathbb{T}^{d}\}$$
is a syndetic subset of $\MNS_{d}(\mathbb{Z})$. 
\end{corollary}
\begin{proof}
Recall $D_{\varepsilon} = \{Y \in \mathcal{K}(\mathbb{T}^{d}) : d(Y,\mathbb{T}^{d}) < \varepsilon\}$, which is open in $\mathcal{K}(\mathbb{T}^{d})$. Let $Y \in \mathcal{I}_{p}$. By Proposition~\ref{prop:randomepsilondense}, there exists $g \in \Gamma$ such that $gY \in D_{\varepsilon}$. Thus
$$\mathcal{I}_{p} \subset \bigcup_{g \in \Gamma} g^{-1}(D_{\varepsilon}).$$
Since $\mathcal{I}_{p}$ is compact, there exists a finite set $F \subset \Gamma$ such that
$$\mathcal{I}_{p} \subset \bigcup_{g \in F}g^{-1}(D_{\varepsilon}).$$
Let $A \in \MNS_{d}(\mathbb{Z})$. Since $T_{A}$ commutes with $T_{p}$ and is finite to one (since $A$ is nonsingular), we have $T_{A}Y \in \mathcal{I}_{p}$. Thus there exists $g \in F$ such that $gT_{A}Y \in D_{\varepsilon}$, and hence $gA \in J_{Y}$. It follows that $J_{Y}$ is syndetic.
\end{proof}

We use this to prove the higher dimensional analog of Mahler's theorem. 
\begin{proof}[Proof of Theorem~\ref{thm:higherdimmahlerintro}]
Let $\alpha \in \mathbb{T}^{d}$ have at least one irrational coordinate and let $\mathcal{M}(\alpha)$ be defined by 
$$\{A \in \MNS_{d}(\mathbb{Z}) : \textrm{ for every } w \in \{0,\ldots,b-1\}^{d}\textrm{, } c_{r}(T_{A}(\alpha)) = w \textrm{ for infinitely many } r\}.$$
Let $\mathcal{B}$ be the set of accumulation points of the set
$$\{T_{b}^{k}(\alpha) : k \ge 1\}.$$
Then $\mathcal{B}$ is $T_{b}$-invariant and infinite, by an argument analogous to the one given in the first paragraph of the proof of Theorem~\ref{th:mahlergeneral}.

Let $\varepsilon = \frac{1}{2b}$. Since $\mathcal{B}$ is also $T_{b}$-invariant, by Corollary~\ref{cor:matricessyndetic} the set
$$J_{\mathcal{B}} = \{A \in \MNS_{d}(\mathbb{Z}) : T_{A}\mathcal{B} \textrm{ is } \frac{1}{2b}\textrm{-dense in } \mathbb{T}^{d}\}$$
is syndetic in $\MNS_{d}(\mathbb{Z})$. It suffices then to show that $J_{\mathcal{B}} \subset \mathcal{M}(\alpha)$. To check this, suppose $A \in J_{\mathcal{B}}$, and so $T_{A}\mathcal{B}$ is $\frac{1}{2b}$-dense in $\mathbb{T}^{d}$. Since $T_{A}$ commutes with $T_{b}$, the set $T_{A}\mathcal{B}$ is the set of accumulation points of $\{T_{b}^{k}T_{A}(\alpha) : k \ge 1\}$. It follows that for every $w=(a_{1},\ldots,a_{d}) \in \{0,\ldots,b-1\}^{d}$, the $T_{b}$-orbit of $T_{A}(\alpha)$ enters the cube $Q_{w}=\prod_{i=1}^{d}(\frac{a_{i}}{b},\frac{a_{i}+1}{b})$ infinitely many times. This implies $c_{r}(T_{A}\alpha)=w$ for infinitely many $r$,  and so $A \in \mathcal{M}(\alpha)$, and hence $J_{\mathcal{B}} \subset \mathcal{M}(\alpha)$.
\end{proof}

\end{document}